# Scientist Mothers Working in Childcare Deserts Are Less Likely to Attain Tenure

Authors: Jun Li[1]*

**Affiliations:**

[1]University of Michigan, Stephen Ross School of Business; Ann Arbor, Michigan 48109, USA.

*Corresponding author. Email: junwli@umich.edu

**Abstract:** Childcare shortages have limited families' access to quality care, with consequences for children's development and parents' careers. We link two decades of childcare supply data to longitudinal surveys of U.S. scientists. We find workplace childcare deserts are widespread: more than half of scientists at universities and research institutions work in deserts, where supply falls significantly short of demand. We further find mothers working in deserts are less likely to attain tenure. They are also more likely to move to industry and part-time positions. Reduced work time partly accounts for the tenure association but other reasons remain. Expanding childcare availability near workplaces may help retain scientific talent.

The US childcare system is characterized by limited access, high costs, and uneven quality. While the COVID-19 pandemic brought the childcare crisis into sharp focus (*1-4*), the problem is decades old (*5, 6*). For years, the limited availability of affordable, high-quality childcare has strained American families (*7*). Yet public investment in childcare remains low relative to that of other high-income countries (*8, 9*). Only a handful of US states provide publicly funded preschool, and such programs primarily serve 4-year-olds only (*10*). For younger children, public care provision is much more limited (*11*), leaving families largely dependent on private providers. The labor-intensive nature of childcare, regulatory requirements, and lack of economies of scale constrain providers' ability to expand capacity. As a result, childcare shortages are widespread (*12, 13*). At the same time, childcare providers often struggle financially despite low wages in the sector and high costs borne by families.

Among the many challenges facing the US childcare system, limited availability is particularly consequential because families cannot use care that is simply unavailable. Before the pandemic, 73% of center-based childcare providers reported turning families away because they were full or maintaining a waiting list in 2019 (*14*). More recent national data also show persistent access difficulties (*15*). Limited access to childcare can reduce labor-force participation (*9, 16, 17*). Moreover, because responsibility for the care of young children continues to fall disproportionately on women (*18, 19*), childcare shortages may also contribute to gender disparities in education and career advancement. These constraints can be especially consequential in scientific careers. A growing literature documents associations between parenthood and researchers' career trajectories and outcomes (*20-23*), yet much less is known about how the environments in which scientists work can support—or constrain—parents' careers. Although parental leave can provide important support around childbirth (*9*), the demands of caregiving extend far beyond the first months of a child's life. Because broadly available public education generally does not begin until kindergarten, families remain heavily dependent on childcare throughout children's early years. In the absence of reliable care, researchers may leave or change professions, reduce their work activities, or reconsider family formation.

Here, we examine workplace childcare availability and its association with scientists' career advancement. In many high-skilled professions, women remain underrepresented at entry and even more so at the top (*24-27*), while disadvantages accumulate over the life course (*28*). Scientific careers at research-intensive universities and institutions require years of doctoral and, sometimes, postdoctoral training, representing substantial investments by both individuals and the institutions and public agencies that support their training. The early years of a scientific career are particularly demanding due to the steep promotion track during this time—tenure track with an almost "up-or-out" trajectory (*29*). These years also frequently overlap with the onset of parenthood. The availability of childcare may therefore be particularly important for whether scientists—especially parents of young children—can sustain the time and effort required for career advancement. Tenure provides long-term employment security and greater research freedom while helping institutions retain established scholars. Yet gender disparities in tenured positions persist despite narrowing gaps in advanced science and engineering education (*30*). We ask two questions. First, how available is childcare near the workplaces of US scientists, and do men and women differ in their exposure to workplace childcare deserts? Second, is working in a childcare desert associated with scientists' likelihood of attaining tenure, and does this relationship differ by gender?

## Results

### *Over Half of US Census Tracts Are Workplace Childcare Deserts*

Childcare provision is fragmented, and licensing records are held by separate state agencies, making national longitudinal data difficult to assemble. We leverage the Infogroup (now Data Axle) US Historical Business Data, a nationwide establishment-level database dating to 1997.

We measure local childcare access at the census-tract level. This does not suggest families cannot access care outside the tract; it merely focuses the analyses on closer care. While our main analyses focus on the workplace tract, robustness analyses also account for childcare availability in neighboring tracts (table S4). A tract is classified as a desert when supply falls below one-third of estimated demand, following Center of American Progress (CAP)'s threshold; Canadian childcare deserts use a similar threshold, i.e., three out of ten (*12, 31*). Unlike prior studies, our focus is workplace rather than residential access, recognizing that parents consider childcare locations jointly with work schedules and commuting routes (*32*).

Workplace childcare deserts remained widespread throughout 1997–2019, with an average of 53.2% (SD = 0.03) of census tracts classified as deserts each year. Fig 1 illustrates their geographic distribution in 2019. CAP's 2018 study found that 51% of Americans lived in residential childcare deserts using licensing data from 22 states (*12*); the measures are different, but both estimates indicate widespread shortages.

### *Long Waiting Lists for Childcare at Leading US Research Institutions*

In 2024, we collected childcare waitlist information at the 50 leading US research institutions, ranked by 2022 research and development expenditures reported by the National Science Foundation. Sixty percent had university-managed centers; most others offered affiliated centers, and two had neither. Among 157 centers, 46.5% reported average waits of one to two years, 27.4% two to three years, and 12.7% six months to one year; one reported a wait of at least three years. Figs. S1 to S3 detail center types, enrollment policies, and waiting times. Although limited to 50 institutions, these data show substantial constraints even at the nation's best-resourced institutions.

### *Childcare Desert Exposure Is Widespread and Same Across Genders, but Tenure Attainment Differs*

Obtained from the Survey of Doctorate Recipient (SDR) under a restricted-use data license, our sample includes 12,629 science, engineering, and health (SEH) doctoral recipients who graduated between 1997 and 2014 and subsequently worked at US universities or research institutions. Individuals are followed up every two to three years, with the most recent wave from 2023.

The sample is nearly evenly divided between men (6,243) and women (6,386). Table S1 reports demographic and professional characteristics. More than half of scientists in our sample were exposed to workplace childcare deserts. Exposure was similar for men and women (57.5% versus 58.0%, $P = 0.592$), indicating that childcare scarcity is widespread but not differentially distributed by gender. Exposure also does not differ by parenthood status or across gender-by-parenthood groups (table S2).

Tenure attainment, however, differs by gender. Overall, 35.5% of scientists attained tenure during the observation period. Women were 5.8 percentage points less likely than men to attain tenure (32.6% versus 38.4%, $P < 0.001$).

Men and women also differ along several other dimensions. Women scientists were less likely to be married or partnered (72.7% versus 76.1%, $P < 0.001$) or have children (39.4% versus 46.3%, $P < 0.001$) below 6 years of age within the first six years upon doctoral graduation. Among scientists who were married or partnered, women were substantially more likely than men to have a spouse or partner working full time (83.5% versus 52.3%, $P < 0.001$). They also differ in age at PhD completion, field of study, characteristics of their PhD and employing institutions, and regional socioeconomic conditions. Our multivariate analyses adjust for these covariates.

### *Childcare Deserts Are Associated with Lower Tenure Attainment among Parents, Particularly Mothers*

We next examine tenure attainment using Cox proportional-hazards models. The models control for all individual, employer, and regional socioeconomic variables. Standard errors are clustered by state. Fig. 2 first presents model-free, cumulative tenure curves (Panel A and B) and then the Cox model estimates (Panel C).

Panel A compares tenure attainment across gender and parenthood groups. Fathers have the highest cumulative tenure attainment over most of the follow-up period, followed by men without children, women without children, and mothers.

Panel B focuses on parents and further separates scientists by childcare-desert exposure. Over most of the follow-up period, fathers working outside childcare deserts have the highest cumulative tenure attainment, followed by fathers working in deserts, mothers working outside deserts, and mothers working in deserts.

Panel C reports estimates from the Cox survival model after accounting for all covariates mentioned above. Results in the first two columns show that women have a lower hazard ratio of attaining tenure, and parenthood further widens the gender difference. Women without children have a 9.1% lower hazard ratio of attaining tenure than men (HR=0.909, SE=0.038, $P = 0.020$). Parenthood is not significantly associated with the tenure-attainment hazard for men (HR=1.031, SE=0.054, $P = 0.562$). The women-and-parenthood interaction is 0.867 (SE=0.059, $P = 0.035$), indicating that the parenthood-associated hazard ratio is 13.3% lower for women than for men.

We now turn to results for childcare desert exposure among scientist parents (Panel C Columns 3 and 4). For fathers, working in a childcare desert is not associated with the likelihood of attaining tenure (HR=1.047, SE=0.087, $P = 0.579$). The women-and-childcare-desert interaction is 0.728 (SE=0.083, $P = 0.005$), indicating that the desert-associated hazard ratio is 27.2% lower for mothers than for fathers. In other words, the association between childcare scarcity and tenure progression is particularly strong among mothers.

Ordinary least squares (OLS) estimates show a similar pattern (table S3). Among scientist parents, women's probability of attaining tenure relative to men's is 6.9 percentage points (pp) lower in childcare deserts than outside them (SE=0.032, $P = 0.034$) during the observation window. The association between desert exposure and tenure attainment among fathers is not statistically significant.

*Childcare Deserts Are Not Associated with Tenure Attainment among Scientists Without Children*

As a placebo test, we estimate the same survival model among scientists without children. We do not observe a statistically significant tenure gap between men and women. Neither do we observe a statistically significant association between childcare-desert exposure and tenure attainment for men or women as shown in Fig. 2C columns 5 and 6.

*Other Services Show No Clear Link to the Tenure and the Gender Gap*

A second placebo test is to examine whether there is an association between labor availability in other service industries such as health, retail, and restaurants, not just childcare, and scientists' tenure attainment. If the association we observe is really attributable to childcare availability, rather than unmeasured, confounding local economic factors, then we should not expect to see an association between tenure rate with labor availability in other service industries. Results shown in Fig. 3 confirms so. Here labor availability is similarly defined across sectors as the percentage of employees working in the corresponding sector out of total number of employees working in a census tract. We see a strong positive association between childcare labor ratio with women scientists' tenure rate, but no association between labor ratios in other sectors and tenure outcome, with the exception that more restaurant labor is associated with lower tenure attainment for men.

***Leakage from Academic Career Paths***

We next examine where researchers go when they do not remain on the path toward tenure. We analyze employment status in the subsequent SDR survey wave and its association with gender and childcare-desert exposure (Fig. 4).

The women-and-childcare-desert interaction indicates that desert exposure has a more positive association with subsequent industry employment for women than for men (4.9 pp, SE=2.1 pp, $P = 0.020$), and a more negative association with employment at universities or research institutions ($-4.8$ pp, SE=2.7 pp, $P = 0.076$). The corresponding interaction for unemployment is small and insignificant ($-0.1$ pp, SE=1.1 pp, $P = 0.913$). We also examine tenure-track appointments and tenure-track or postdoctoral appointments at universities or research institutions. For tenure-track appointments, the interaction is negative but does not reach conventional significance levels ($-5.2$ pp, SE=3.5 pp, $P = 0.143$). The interaction is positive for subsequent part-time employment overall (5.1 pp, SE=3.0 pp, $P = 0.094$) and for part-time employment at universities or research institutions (4.1 pp, SE=2.0 pp, $P = 0.042$).

Together, these patterns suggest that childcare shortages are associated with movement away from tenure-track academic careers and toward industry as well as part time positions with no statistically significant association with unemployment.

***Lack of Childcare Means Less Time Available for Mothers, while Other Mechanisms Remain***

Why is limited childcare associated with women's tenure attainment more than men's? Pandemic survey evidence links caregiving for young children to larger reductions in scientists' research time, especially for mothers (*1*). One possibility, therefore, is that childcare shortages reduce the time mothers can devote to work.

We test this mechanism using the SDR question, "During a typical week on your principal job, how many hours did you work?" The average reported workweek in our sample is 48.7 hours (SD=12.5 hours).

Parenthood is associated with approximately 1.6% fewer working hours for men (SE=1.0%, $P = 0.102$), albeit not reaching conventional significance levels (Table 1, column 1). For women, parenthood is associated with an additional reduction of 13.5% in working hours relative to men (SE=1.3%, $P < 0.001$). Scaling these estimates by the sample mean gives 0.8 fewer hours per week for men and an additional 6.6 fewer hours for women than men.

Among parents, working in childcare deserts is not significantly associated with fathers' working hours (0.009, SE=0.011, $P = 0.403$). For mothers, desert exposure is associated with an additional reduction of 4.9% in working hours relative to fathers (SE=2.3%, $P = 0.038$; Table 1, column 2). At the sample mean, it corresponds to 2.4 fewer hours for mothers than fathers.

Reduced work time explains only part of the tenure disparity. When controlling for working hours, women scientists working in childcare deserts are still less likely to attain tenure relative to men. Other mechanisms may therefore matter. Parents facing unreliable childcare may shift work to less desirable hours, experience sleep loss and fatigue (*33, 34*), or face frequent work interruptions. Multitasking and reduced synchronization with collaborators may further lower productivity even when total working hours change little. Because working hours are self-reported and capture total hours rather than their timing or fragmentation, our data cannot distinguish among these mechanisms.

**Discussion**

About 46,000 people earn science, engineering, and health doctorates from US institutions annually (*35*), but relatively few attain tenure. Although women constitute more than half of tenure-track assistant professors and non-tenure-track professors, they hold only 45% of associate and 32.5% of full professorships at research universities and institutions (*36*). While recent literature shows that parenthood is associated with researchers' career outcomes (*21, 22*), the implication is not that researchers, and American families, should stop having children. Rather, society should find ways to better support young families.

Workplace childcare-desert exposure is widespread and associated with lower tenure attainment among scientist parents, particularly mothers, after adjustment for individual, family, employment, and local socioeconomic characteristics and state, graduation-year, and response-year fixed effects. The findings are robust to alternative specification under OLS (table S3) and extending deserts measure to neighboring tracts (table S4). We observe the association among parents but not scientists without children, and other local services show no comparable gender-specific pattern. We also find no evidence that workplace childcare-desert exposure differs systematically by gender, parental status, or the timing of parenthood. Together, these comparisons support a relationship specific to childcare availability rather than simply differences across local labor markets.

Our analyses also provide evidence on a potential mechanism. Part of the gender difference associated with childcare deserts can be explained by reductions in mothers' time available for work. The estimated contribution of work time, however, is moderate, suggesting that additional mechanisms may be important. Limited childcare availability may also contribute to burnout, interruptions to work, greater demands for multitasking, and difficulties coordinating schedules with collaborators. Understanding these mechanisms represents an important direction for future research.

The consequences of limited childcare availability extend well beyond scientific careers. Many US families have faced persistent shortages of childcare, with potential implications for labor-force participation, productivity, and well-being. Policies to expand childcare availability and access include investments in infrastructure and provider capacity, financial support for providers, family subsidies (37), and publicly funded or universal programs (*38*). For example,

Michigan has invested in opening and expanding childcare facilities (*39*), Louisiana provides tax credits to providers and workers (*40*), and the federal Child Care and Development Fund subsidizes childcare for eligible families (*41*). States have also expanded publicly funded early childhood education (*10*); most recently, New Mexico established no-cost universal childcare regardless of family income (*42*). The relative effectiveness of these approaches, and the settings in which they are most effective, remain important research questions.

Employers also have an important role. A National Academies consensus report recommends stronger institutional support for family caregivers in STEMM (*43*). Universities and other employers could consider approaches such as on-site childcare, backup-care programs, partnerships with local providers, or childcare subsidies (*43*). Our findings underscore the importance of considering childcare availability near workplaces when designing institutional support for scientists.

Limitations of this study also point to opportunities for future research. First, our findings are primarily correlational. Although the analyses account for a wide range of potential confounders and are robust to comparisons involving scientists without children and the availability of other local services, unobserved differences across researchers or locations may remain. Future studies that exploit plausibly exogenous changes in childcare availability could provide stronger evidence on the causal effects of childcare access on scientific careers. Second, our childcare data capture center- and group-based care, which represent more than 80% of all licensed slots (*44*). We do not observe family-based care or informal care arrangements. Incorporating these alternatives care options will provide a more comprehensive view of childcare availability and its potential impact.

Finally, childcare availability is only one dimension of the broader childcare challenge facing families. Affordability and quality are also important. Increasing supply without making it affordable will not truly increase accessibility (*45, 46*). Meanwhile, making childcare universal without investing in quality could have negative short- and long-term impacts on children (*16, 47, 48*). A major obstacle to studying these dimensions at scale is the lack of granular, nationwide, longitudinal data on childcare prices and quality. Building such data infrastructure would enable future research to examine availability, affordability, and quality jointly and provide a more complete understanding of the causes and consequences of the US childcare crisis.

**Acknowledgments:**

The author thanks seminar participants at Dartmouth College, Georgia Institute of Technology, Imperial College London, London Business School, Massachusetts Institute of Technology, Pennsylvania State University, Syracuse University, Tilburg University, University College London, University of Miami, University of South Carolina, Western University, and Yale University for their helpful comments and feedback on this research. The author also thanks Elizabeth Hwang, Claire Miech, Maddison Linker, and Liam Belluche for their assistance with data collection.

## Funding:

No funding was received for this research.

**Author contributions:** The paper is single authored. The author confirms sole responsibility for the following: study conception and design, data collection, analysis and interpretation of results, and manuscript preparation.

**Conceptualization:** JL

**Methodology:** JL

**Investigation:** JL

**Formal analysis:** JL

**Writing – original draft:** JL

**Competing interests:**

The author declares no competing interests.

**Data, code, and materials availability:** The study utilizes the following datasets:

(1) The Infogroup US Historical Business Data (1997-2020) are accessed through University of Michigan's License with Infogroup (now DataAxle). The data can be accessed from University of Michigan's Library:

https://search.lib.umich.edu/databases/record/9111

(2) Survey of Doctorate Recipients (SDR) Linked with the Doctorate Records File (DRF) data files are accessed through National Center for Science and Engineering Statistics (NCSES)'s restricted-use data license. The data are stored and analyzed within NCSES's Virtual Data Enclave. Researchers should follow the instructions listed on NCSES's webpage below to obtain the license and to access restricted use data:

https://ncses.nsf.gov/about/licensing

(3) The Integrated Postsecondary Education Data System (IPEDS) data are publicly available and can be downloaded directly from IPEDS' webpage below:

https://nces.ed.gov/ipeds/use-the-data

(4) Census Bureau Data are publicly available and can be downloaded directly from the Census's webpages below:

https://www.census.gov/programs-surveys/popest.html

https://www.census.gov/programs-surveys/cbp/data/datasets.html

https://www.census.gov/library/publications/2011/compendia/usa-counties-2011.html

(5) Presidential Elections in US Counties (1868–2020) data are publicly available through the authors' Harvard Dataverse repository:

https://doi.org/10.7910/DVN/DGUMFI

The final merged data are stored and analyzed within the NCSES's Virtual Data Enclave per NCSES's licensing policy.

*Code availability.* All STATA code (version 18) are available within the NCSES's Virtual Data Enclave per NCSES's licensing policy.

**Supplementary Materials**

Materials and Methods

Figs. S1 to S3

Tables S1 to S4

References (49–56)

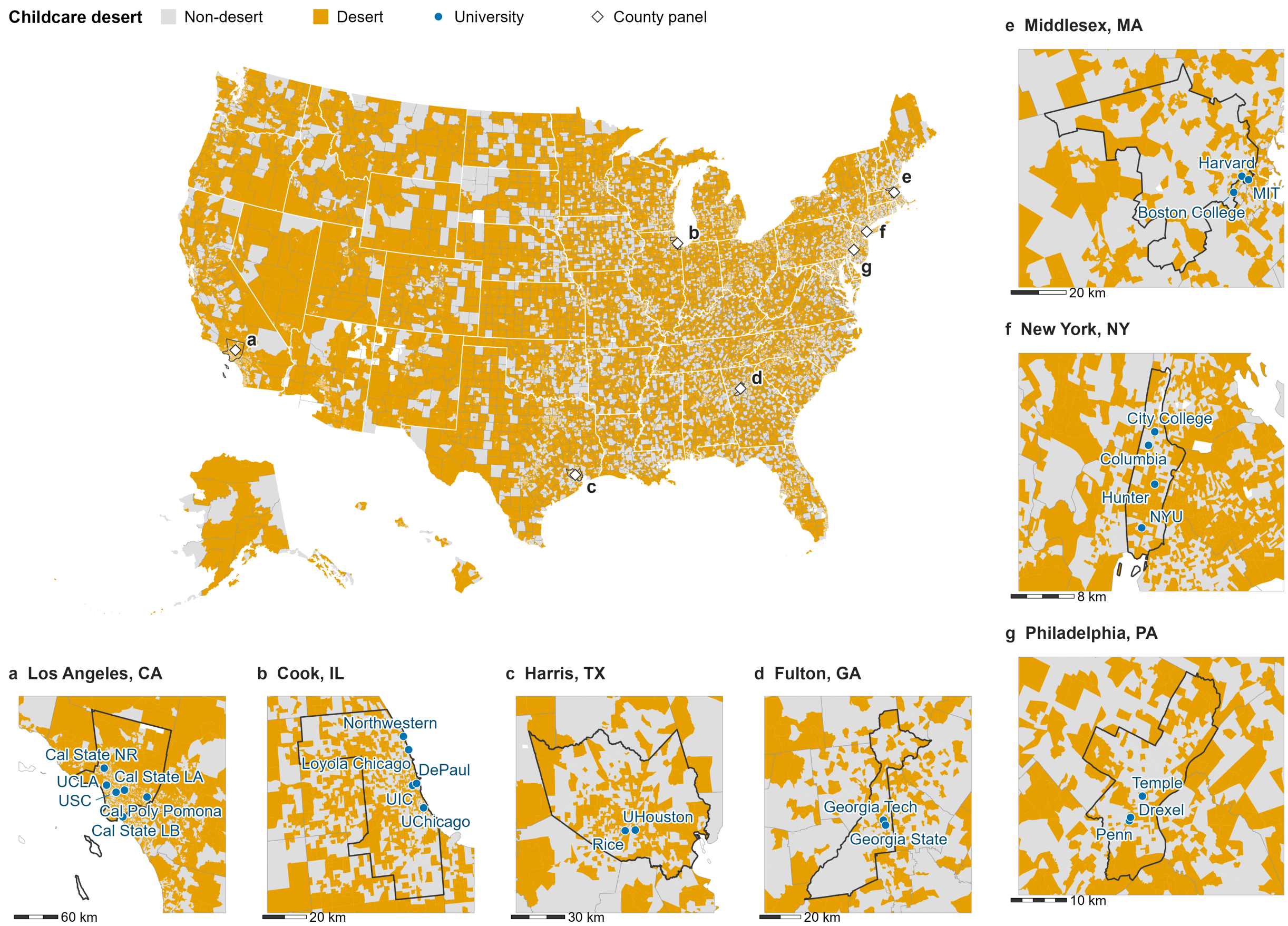

**Fig. 1. Workplace childcare deserts at the census-tract level, 2019.** Workplace childcare deserts (shown in orange) are defined as census tracts in which there is no or very low childcare supply, i.e., a supply-to-demand ratio less than one third. Workplace childcare supply is approximated as proportional to the number of employees working in childcare facilities (primary SIC code: 8351) in a census tract based on the Infogroup (now Data Axle) US Historical Business Data. Workplace childcare demand is approximated as proportional to the total number of employees working in a census tract. Refer to the Methods section for a detailed definition.

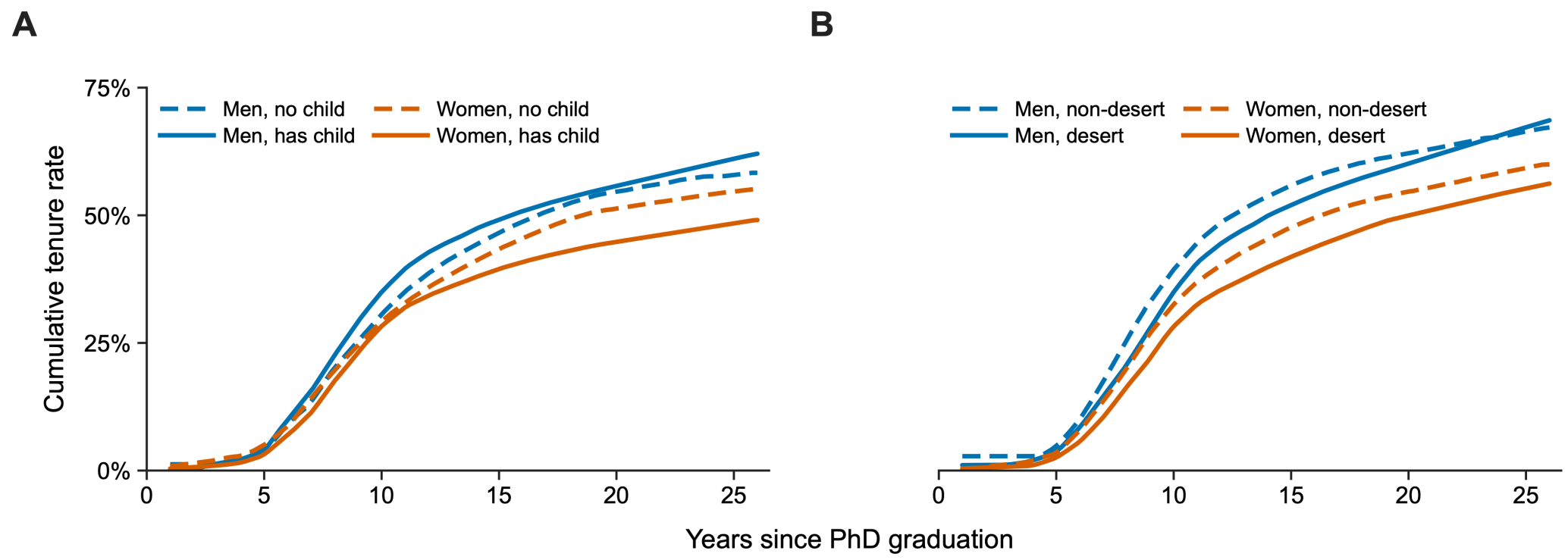

**C Cox proportional-hazards estimates**

| | All scientists | | Parents | | Nonparents | |
|---|---|---|---|---|---|---|
| | (1) | (2) | (3) | (4) | (5) | (6) |
| Women | 0.856*** (0.034) | 0.909** (0.038) | 0.795*** (0.059) | 0.962 (0.081) | 0.934 (0.050) | 0.922 (0.087) |
| Parenthood | | 1.031 (0.054) | | | | |
| Women × parenthood | | 0.867** (0.059) | | | | |
| Childcare desert | | | | 1.047 (0.087) | | 0.982 (0.086) |
| Women × childcare desert | | | | 0.728*** (0.083) | | 1.015 (0.125) |
| Observations | 53,512 | 53,512 | 21,713 | 21,674 | 24,869 | 24,831 |
| Log pseudolikelihood | −1,147,532.3 | −1,147,472.7 | −473,200.2 | −472,964.0 | −485,799.1 | −485,209.6 |

Hazard ratios; standard errors in parentheses. *** $P < 0.01$; ** $P < 0.05$; * $P < 0.10$.

**Fig. 2. Gender, parenthood, and childcare deserts in tenure promotion. (A)** Cumulative tenure rates by gender and parenthood. **(B)** Cumulative tenure rates by gender and workplace childcare-desert status. **(C)** Hazard ratios from six specifications of Cox proportional hazards models estimated for all scientists, parents, and nonparents. Hazard ratios below 1 indicate a lower instantaneous rate of tenure promotion. All models include the full set of individual, employer, and regional controls listed in table S1. Individual controls include race and ethnicity, citizenship, marital or partnership status, indicators for spousal full-time and part-time employment, age at PhD completion, whether the doctoral institution is classified as having the highest research activity (Carnegie Classification of Institutions of Higher Education), PhD field of study (NSF classification), year of PhD completion, and year of first SDR response. Employer controls include whether the employing institution is classified as having the highest research activity (Carnegie Classification of Institutions of Higher Education); log faculty counts overall and at the full, associate, and assistant professor ranks; percentages of women faculty overall and by rank; and log undergraduate and graduate student counts. Regional controls include commuting-zone measures of log population density, log business density, log employment density, log annual payroll density, log median household income, the unemployment rate, and log population density of children aged 0–5, together with the county Democratic presidential vote share. All models also include state fixed effects. Parenthood is defined as having at least one child younger than 6 during the first six years after PhD completion. The parent subsample uses this definition. The nonparent subsample in columns (5) and (6) comprises scientists with no children younger than 18 for clean comparison. Standard errors, clustered at the state level, are reported in parentheses. Observation counts are follow-up records rather than unique respondents.

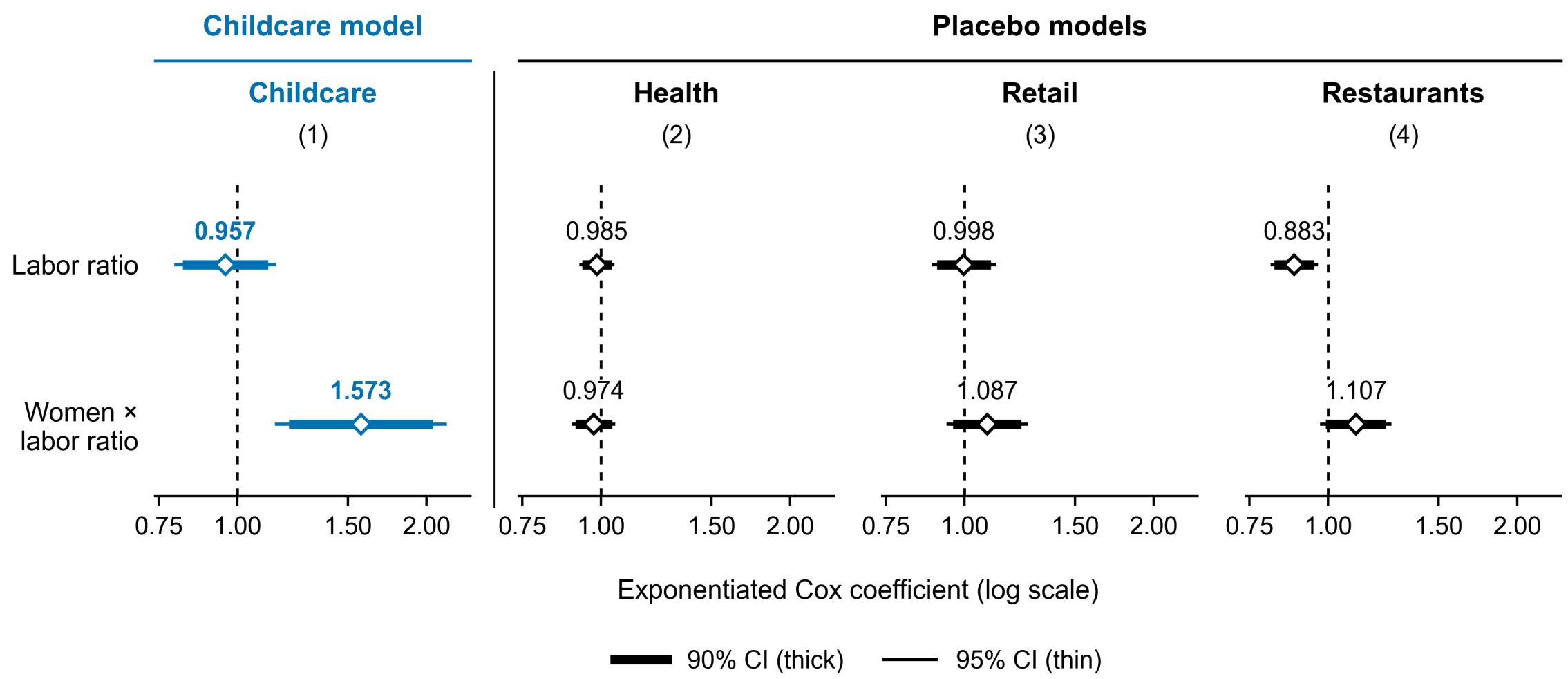

**Fig. 3. Placebo tests using labor ratios in other industries.** Hazard-ratio plots from four Cox proportional hazards models of tenure promotion are shown: childcare (1), health (2), retail (3), and restaurants (4). Each model includes the women indicator, the industry labor ratio and its interaction with women, as well as all control variables as previously listed in Fig. 2C. Only labor ratio and its interaction with the women indicator are shown in the figure. Industry labor ratio is defined as the total number of employees in that specific industry divided by the total number of employees in a census tract and is log scaled. Diamonds mark point estimates; thick and thin segments denote 90% and 95% confidence intervals, respectively. The dashed line marks a ratio of 1. Blue identifies the childcare model; black identifies the three placebo models. All models are displayed on identical logarithmic scales.

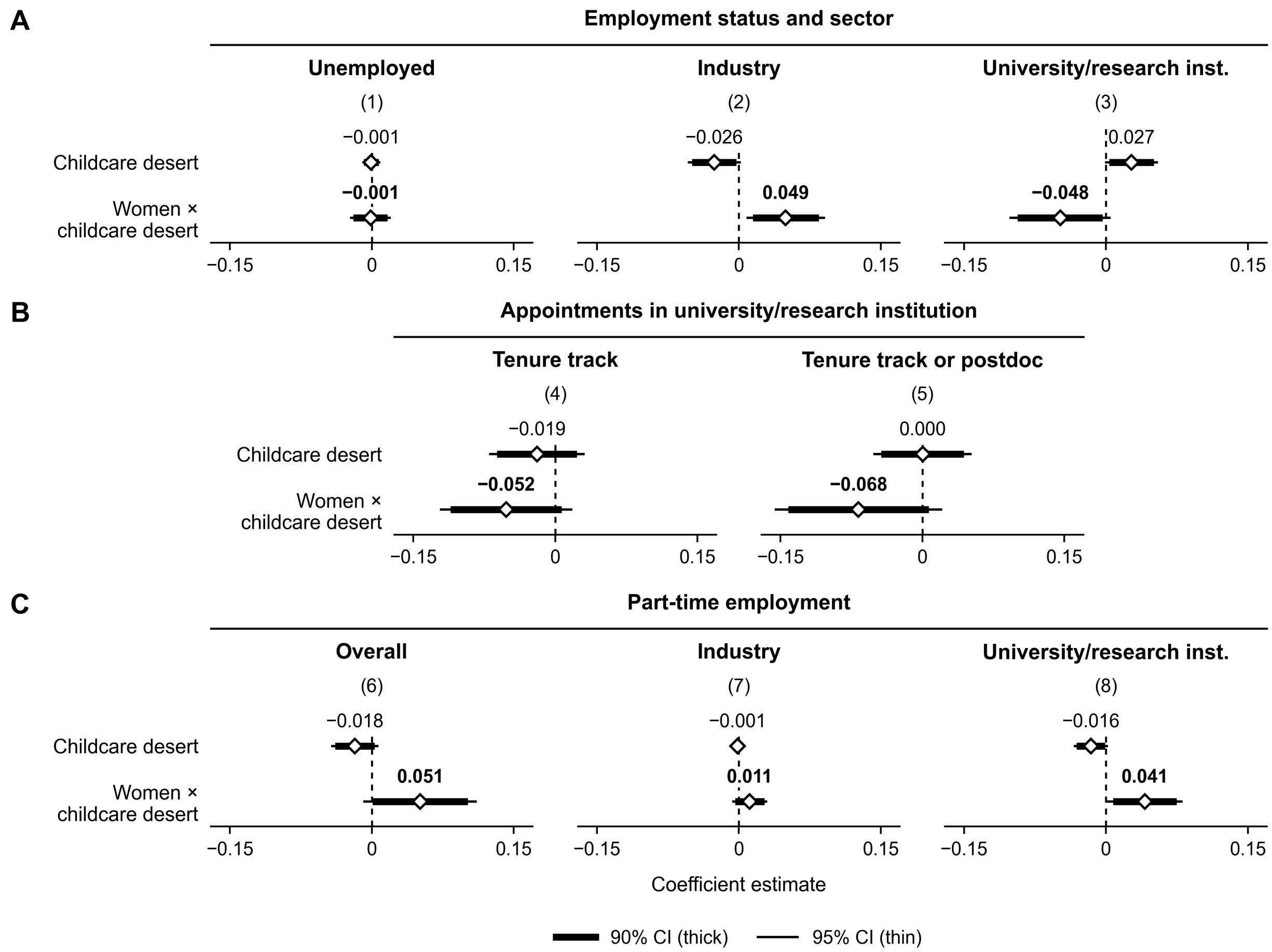

**Fig. 4. Childcare deserts and subsequent employment outcomes.** Results from OLS regressions. The outcome is employment status in the subsequent SDR wave. **(A)** Broad employment status and sector outcomes: unemployment, employment in industry, and employment at universities or research institutions. **(B)** Appointment outcomes within universities or research institutions: tenure track, and tenure track or postdoctoral appointment. **(C)** Part-time employment: overall, in industry, and at universities or research institutions. Each model includes the women indicator, a workplace childcare-desert indicator and its interaction with women, as well as all control variables as previously listed in Fig. 2C. Only coefficients for the childcare-desert indicator and its interaction with women are displayed. Diamonds mark point estimates; thick and thin segments denote 90% and 95% confidence intervals, respectively. The dashed line marks zero.

**Table 1. Gendered associations of parenthood and childcare deserts with weekly working hours and tenure promotion.**

| Sample | All scientists | Parents | Parents |
|---|---|---|---|
| Outcome | Log weekly hours | Log weekly hours | Tenure promotion |
| Model | OLS | OLS | Cox |
| Reported statistic | Coefficient | Coefficient | Hazard ratio |
| | (1) | (2) | (3) |
| Women | $-0.021^{*}$ | $-0.123^{***}$ | 1.072 |
| | (0.011) | (0.012) | (0.085) |
| Parenthood | $-0.016$ | | |
| | (0.010) | | |
| Women × parenthood | $-0.135^{***}$ | | |
| | (0.013) | | |
| Childcare desert | | 0.009 | 1.035 |
| | | (0.011) | (0.082) |
| Women × childcare desert | | $-0.049^{**}$ | $0.743^{***}$ |
| | | (0.023) | (0.074) |
| Log weekly working hours | | | $3.362^{***}$ |
| | | | (0.383) |
| Observations | 12,150 | 5,011 | 21,674 |
| Adjusted $R^2$ | 0.088 | 0.110 | |
| Log pseudolikelihood | | | $-470{,}309.2$ |

*Notes.* Columns (1) and (2) report OLS coefficients for the natural logarithm of reported weekly working hours. Column (3) reports hazard ratios from a survival model of tenure promotion, controlling for the logarithm of average weekly working hours. All models include the same controls described in Fig. 2C. Standard errors are in parentheses and are clustered at the state level.

$^{***}P < 0.01$; $^{**}P < 0.05$; $^{*}P < 0.10$.

**Supplementary Materials for**

# Scientist Mothers Working in Childcare Deserts Are Less Likely to Attain Tenure

Jun Li

*Corresponding author. Email: junwli@umich.edu

This file includes Materials and Methods, Figs. S1 to S3, and Tables S1 to S4. References are numbered in the single References and Notes list in the main manuscript. References 49–56 are cited only in these Supplementary Materials.

## Materials and Methods

### *Ethics*

The University of Michigan's Institutional Review Board approved the study as Not Regulated (category: Research Involving Coded Private Information/Research, and Research Involving De-identified Biological Specimens or Information). The PI (and the Regents of the University of Michigan) obtained the license for restricted-use data from NCSES in February 2022 following the NCSES's licensing procedure. The license is renewed in June 2026 and is effective till June 2029. The PI is the only person with approved access to the licensed data. The data are stored, accessed and analyzed within the NCSES Virtual Data Enclave. The PI has completed annual data security training, undergone periodic data security audits, and maintained continued compliance per the license agreement. The tables and figures presented in this manuscript comply with the data privacy and confidentiality guidelines of the NCSES and have been reviewed and approved by the NCSES Data Enclave staff for export and publication.

### *Data*

This study uses the following datasets:

(1) *Infogroup (now Data Axle) US Historical Business Data*: This dataset contains detailed facility-level information on US businesses and organizations since 1997. It includes data on the facility's name, location, its parent company, industry classifications (with primary SIC and NAICS codes, and up to six secondary SIC codes), number of employees, and sales volume. This dataset is integral to our research because it enables us to identify childcare facilities, as well as their approximate capacities, using the primary SIC code 8351. Childcare facilities under this classification include childcare centers, day care centers, group day care centers, Head Start centers, nursery schools, and preschool centers.

To ensure data quality, we compare the total number of childcare facilities and employees identified from the Infogroup database with those reported in the US Census Bureau's County Business Patterns (CBP) database (available only at the county level). Across all years in this study, the numbers of childcare facilities and employees in the two databases are highly correlated, with correlations of 96% to 97%.

(2) *Survey of Doctorate Recipients (SDR) Linked with the Doctorate Records File (DRF)*: The SDR, a biennial survey, targets individuals under 76 years old who have earned a research doctoral degree in a science, engineering, or health field from a US institution. The survey's population includes roughly 1.2 million individuals, with about 120,000 participants in the sample. It gathers data on demographics, education, family, and career histories, and is co-sponsored by the National Center for Science and Engineering Statistics (NCSES) within the National Science Foundation and the National Institutes of Health. The most recent data available are from 2023.

The DRF complements this by being a cumulative record of all US research doctorate recipients since the 1920s. Updated annually on the basis of the Survey of Earned Doctorates (SED), DRF provides data on educational history, demographic characteristics, and post-graduation plans. The SED is sponsored by the NCSES, as well as three federal agencies (the National Institutes of Health, US Department of Education, and the National Endowment for the Humanities), and the most recent data available are from the academic year 2023. Access to these data is obtained through an NSF/NCSES restricted-use data license.

(3) *Integrated Postsecondary Education Data System (IPEDS)* Data: IPEDS is a system that gathers data from every US college and university through twelve interrelated survey components conducted annually. In this study, we use institutional characteristics, enrollment statistics, and human resources metrics from the database. Each postsecondary institution in the survey is designated a unique UNITID, which is used to link the IPEDS data and the SDR data based on employer institution UNITIDs.

(4) *Census Bureau Data*: The Census Bureau provides annual estimates of population, economic establishments, total employment, median household income, and unemployment rate at the county level, as well as the land area of each county.

(5) *Presidential Elections in US Counties (1868–2020)*: The dataset published in *Electoral Studies* allow us to obtain county-level Democratic presidential vote shares (*49*).

***Data Merging Process***

We take the following three steps to integrate the datasets mentioned above.

Step 1: Our initial step involves combining the IPEDS data, which include geo-coded institutional locations, with census tract level childcare availability data derived from the Infogroup US Historical Business data. IPEDS data provide geo-coded locations for institutions. Using R's tigris package, we first identify the 2020 census tracts associated with these geo-locations. When necessary, we convert the 2020 census tracts to 2010 and 2000 census tracts using the relationship files provided by the Census Bureau. Then, we merge the IPEDS data with the childcare availability data based on each institution's census tract location.

Step 2: In the second step, we proceed to link the combined IPEDS and childcare data with the regional economic metrics at the commuting zone level.

Step 3: For the final step, we integrate the previously merged datasets with the SDR/DRF dataset. The merging is conducted at the institution level using unique institutional identifiers, known as UNITIDs.

***Sample***

We restrict the SDR sample to individuals who received doctoral degrees between 1997 and 2014 to allow sufficient time to observe tenure outcomes. With most recent available survey data extending through 2023, the 2014 cutoff allows seven years for progression to tenure, through 2021, plus two additional years of potential follow-up to accommodate pandemic-related delays. We retain respondents with at least one survey response within six years of PhD completion and a final observed response at least of seven years after completion. This selection criterion enables us to track employment locations shortly after graduation and observe subsequent tenure outcomes. We acknowledge that tenure track lengths vary across institutions as well as due to individual extensions. Therefore, our sample may include recent graduates who are still in tenure-track positions and have not yet reached the tenure evaluation phase as of the 2023 survey.

We further restrict the sample to respondents initially employed at four-year colleges or universities, medical schools, and university-affiliated research institutions in the US and link their employer institutions to IPEDS and childcare-availability data. The final linked sample contains 12,629 respondents. Unmatched individuals either lack valid employment institution codes or are employed by institutions located in US territories (e.g., Guam, Puerto Rico) for which no childcare data are available.

***Empirical Specifications***

Our primary analyses use Cox proportional hazards models to examine associations of gender, early-career parenthood, and workplace childcare deserts with time to tenure attainment (Fig. 2). Let $h_i(t)$ denote respondent $i$'s instantaneous rate of attaining tenure at $t$ years since PhD completion, conditional on not yet attaining tenure, and let $h_0(t)$ denote the unspecified baseline hazard. Tenure attainment is dated at the first SDR response in which the respondent reports being tenured. Respondents without an observed tenure event are right-censored when their observed follow-up ends.

The survival data contain multiple records per respondent, representing successive intervals of observed follow-up. Each respondent contributes records while at risk and can experience the tenure event only once. Parenthood captures whether the respondent reported having a child aged 0–5 within the first six years after PhD completion and is treated as a fixed early-career characteristic. Childcare-desert exposure and the control variables are also held fixed across a respondent's follow-up records. For parents, childcare desert exposure and the associated controls are evaluated using the first SDR response within this early-career window at which the respondent was employed at an academic institution and reported a young child. Thus, the multiple records describe follow-up over time, while the covariates retain their selected early-career values.

Let $F_i$ denote the female indicator, $P_i$ the early-career parenthood indicator, and $D_i$ the workplace childcare-desert indicator. The vector $X_i$ contains the individual, employer institution, and regional socioeconomic controls described in the Variables subsection, as well as indicators for PhD graduation year, first SDR response year, and state fixed effects.

For the full sample, we first estimate a model with the female indicator and controls, then add parenthood and its interaction with gender:

$$h_i(t) = h_0(t)\exp(\beta F_i + X_i'\gamma), \qquad \text{(S1)}$$

$$h_i(t) = h_0(t)\exp(\beta F_i + \delta P_i + \eta F_i P_i + X_i'\gamma). \qquad \text{(S2)}$$

For the parent subsample, we estimate the corresponding gender-only model in Equation S1 and a model that includes childcare-desert exposure and its interaction with gender:

$$h_i(t) = h_0(t)\exp(\beta F_i + \delta D_i + \eta F_i D_i + X_i'\gamma). \qquad \text{(S3)}$$

Parameters and baseline hazards are estimated separately for each specification and sample. We also estimate the same pair of specifications among scientists without children younger than 18 as a placebo analysis. For this nonparent subsample, childcare desert exposure and the associated controls are evaluated using the first SDR response within this early-career window at which the respondent was employed at an academic institution. Standard errors are clustered at the state level.

As a robustness check, we estimate linear probability models by ordinary least squares (OLS), using the same covariates and interaction structure (table S3):

$$y_i = \alpha + \beta F_i + X_i'\gamma + \epsilon_i, \qquad \text{(S4)}$$

$$y_i = \alpha + \beta F_i + \delta P_i + \eta F_i P_i + X_i'\gamma + \epsilon_i, \qquad \text{(S5)}$$

$$y_i = \alpha + \beta F_i + \delta D_i + \eta F_i D_i + X_i'\gamma + \epsilon_i. \qquad \text{(S6)}$$

Here, $y_i$ equals 1 if the respondent ever reported being tenured during observed follow-up, and 0 otherwise. The gender-only (S4) and parenthood (S5) specifications are estimated for the full sample; the gender-only (S4) and childcare-desert (S6) specifications are estimated for the parent and nonparent subsamples. The OLS coefficients describe differences in the probability of observed tenure attainment (*50*) during observation window. Since the follow up window is different for individual graduated in different years, the inclusion of the graduation year fixed effects address this, but the estimate is best interpreted as average differences for the sample and during the sample period. Standard errors are again clustered at the state level.

To study leakage from academic career paths (Fig. 4), we estimate separate linear probability models by OLS for employment outcomes reported in the subsequent SDR survey:

$$y_{i,t+1} = \alpha + \beta F_i + \delta D_{i,t} + \eta F_i D_{i,t} + X_{i,t}'\gamma + \epsilon_{i,t+1}. \qquad \text{(S7)}$$

Here, $t$ indexes SDR survey waves, which are generally two years apart and occasionally three years apart. For each outcome, $y_{i,t+1}$ equals 1 if respondent $i$ reports the specified employment status in wave $t+1$, and 0 otherwise. Only respondent-wave pairs with an observed response in the subsequent survey enter these models. $D_{i,t}$ denotes workplace childcare-desert exposure in wave $t$, $F_i$ is the female indicator, and $X_{i,t}$ contains the same individual, employer institution, and regional socioeconomic controls described above, evaluated in wave $t$.

We estimate this model separately for unemployment; employment in industry or at universities or research institutions; tenure-track appointments and tenure-track or postdoctoral appointments at universities or research institutions; and part-time employment overall, in industry, or at universities or research institutions. The interaction coefficient $\eta$ captures the difference between women and men in the association of desert exposure with the subsequent employment outcome. Standard errors are clustered at the state level.

All Cox and OLS regressions use SDR survey weights adjusted for selection into our study sample. For each respondent, the base weight is the SDR-provided survey weight for the year of the first SDR response, which already accounts for survey nonresponse and changes in the sampling panel over time. We further adjust this weight by multiplying it by the inverse of the estimated probability of selection into the study sample, conditional on gender and the full set of individual-level characteristics described in the Variables subsection.

### *Variables*

We define the outcome and main explanatory variables below, followed by three groups of control variables: individual characteristics, employer institution characteristics, and regional socioeconomic characteristics. Table S1 reports summary statistics by gender and for the full sample.

**Tenure attainment, parenthood, and gender.** Tenure attainment (**Tenured**) equals 1 if the respondent ever reported being tenured in the SDR, and 0 otherwise. Parenthood (**Has Child**) equals 1 if the respondent reported having at least one child aged 0–5 in an SDR response within the first six years after PhD completion, and 0 otherwise. Gender is captured by the **Female** indicator obtained directly from the SDR.

**Workplace childcare deserts.** A workplace childcare desert is a census tract where the childcare supply-to-demand ratio falls below one third. We measure childcare deserts at the census tract level as childcare demand is local. The median land area of a census tract is about 1.8 square miles, and it varies based on local population and economic densities (25th percentile: 0.7 square miles; 75th percentile: 12 square miles). Our use of census tracts aligns with that use in the definition of childcare deserts by the CAP (*12*), and the unit used for food deserts and pharmacy deserts (*51, 52*). Tract boundaries describe closer access rather than limit where families can obtain childcare. Our main analyses use desert exposure in the workplace tract. Spatial robustness models jointly include separate desert indicators for the focal and neighboring tracts and each indicator's interaction with gender (table S4).

We calculate childcare supply in a census tract as the total number of employees working in childcare facilities in the tract. Because staff-to-child ratios are regulated, employment provides a proxy for facility capacity, although these ratios vary across states and over time. Our model account for state fixed effects and time (graduation year and response year) fixed effects. The measure approximates capacity rather than counting licensed slots directly. To quantify childcare demand, we note that workplace locations in our study are primarily university campuses, where the majority of residents are college students. Therefore, the primary source of childcare demand in these areas comes from the workforce rather than the residents. We thus define childcare demand as proportional to the number of employees working in a census tract, again using the Infogroup database. This workplace-based denominator differs from CAP's residential measure, which uses the number of young children living in the tract. Ideally, an accurate measure of demand would require counting the number of children aged 0 to 5 of these employees. Unfortunately, data of this granularity (tract level) and longitudinal span are unavailable. The most relevant is the American Community Survey (ACS) 2000–2005 snapshot which contains a five-year estimate of the total number of workers with children under age 6 at the census tract level. We find the two estimates—our measure derived from the Infogroup database and the ACS five-year estimates—are highly consistent, with correlations as high as 87.0% at the census tract level and 98.9% at the county level.

We then calculate the childcare supply-to-demand ratio for each census tract and adopt a threshold to define childcare deserts. The threshold is chosen such that a childcare desert is characterized as an area with no or extremely low child care supply, i.e., less than one third of demand.

Specifically,

$$\begin{gathered}\text{Workplace Childcare Supply-to-Demand Ratio} \\ \approx \frac{\begin{gathered}\text{Childcare Employee Count} * \text{\%Direct Labor} \\ * \text{Children-to-Staff Ratio}/\text{\%Center Based Childcare Slots}\end{gathered}}{\begin{gathered}\text{Total Employee Count} \\ * \text{Average Number of Children Aged 0-5yrs Per Employee}\end{gathered}} \\ \approx a * \frac{\text{Number of Childcare Employees}}{\text{Total Number of Employees}}.\end{gathered}$$

That is, our workplace childcare supply-to-demand ratio is defined as proportional to the percentage of the workforce in a census tract working in the childcare industry. Intuitively, a greater percentage of the workforce working in the childcare industry indicates higher childcare availability in an area, and a lower percentage indicates lower availability. The multiplier $a$ can be approximated using the following statistics: (1) average percentage of direct labor in the childcare industry (84%)(*53*), (2) average child-to-staff ratio (8:1)(*54*), (3) market share of center or group based childcare slots (85%)(*44*), and (4) average number of children aged 0-5 years per worker (1/9), which is calculated by the average population aged 0-5 years /(average population aged between 20 and 65 years * (1- average unemployment rate)).

We then define workplace childcare deserts as follows:

$$\begin{gathered}\text{Workplace Childcare Desert} \\ = \begin{cases} 1, & \text{if Workplace Childcare Supply-to-Demand Ratio} < \frac{1}{3}, \\ 0, & \text{otherwise.} \end{cases}\end{gathered}$$

**Individual characteristics.** Individual controls include race and ethnicity, U.S. citizenship, marital or partnership status, and indicators for spousal full-time and part-time employment, obtained from the SDR. We also include educational characteristics: age at PhD completion, and PhD field of study, obtained from the DRF, together with an indicator for whether the doctoral institution is classified as having the highest research activity under the Carnegie Classification of Institutions of Higher Education. We also include fixed effects for year of PhD completion and the first SDR response.

**Employer institution characteristics.** Employer controls include an indicator for whether the employing institution is classified as having the highest research activity under the Carnegie classification; faculty counts overall and at the full, associate, and assistant professor ranks; percentages of women faculty overall and at each rank; and undergraduate and graduate student counts. Faculty and student counts are log transformed. These institutional characteristics are obtained from IPEDS.

**Regional socioeconomic characteristics.** Regional controls include commuting-zone densities of population, business establishments, employment, annual payroll, and children aged 0–5, together with median household income and the unemployment rate. Density measures are defined as the corresponding counts or totals divided by land area. Density measures and median household income are log transformed. We also include the county Democratic presidential vote share as a regional control and state fixed effects. The Democratic presidential vote share is calculated as the number of votes cast for the Democratic presidential nominee divided by the total number of votes cast in the county in the most recent presidential election, multiplied by 100.

We aggregate county-level socioeconomic measures to commuting zones to better capture characteristics of the local labor markets. County boundaries often reflect political boundaries, whereas commuting zones are intended to reflect the local economies where people live and work (*55, 56*). Delineations of commuting zones before 2010 are obtained from USDA's Economic Research Service website, and those after 2010 are obtained from the Pennsylvania State University's Labor-sheds for Regional Analysis website.

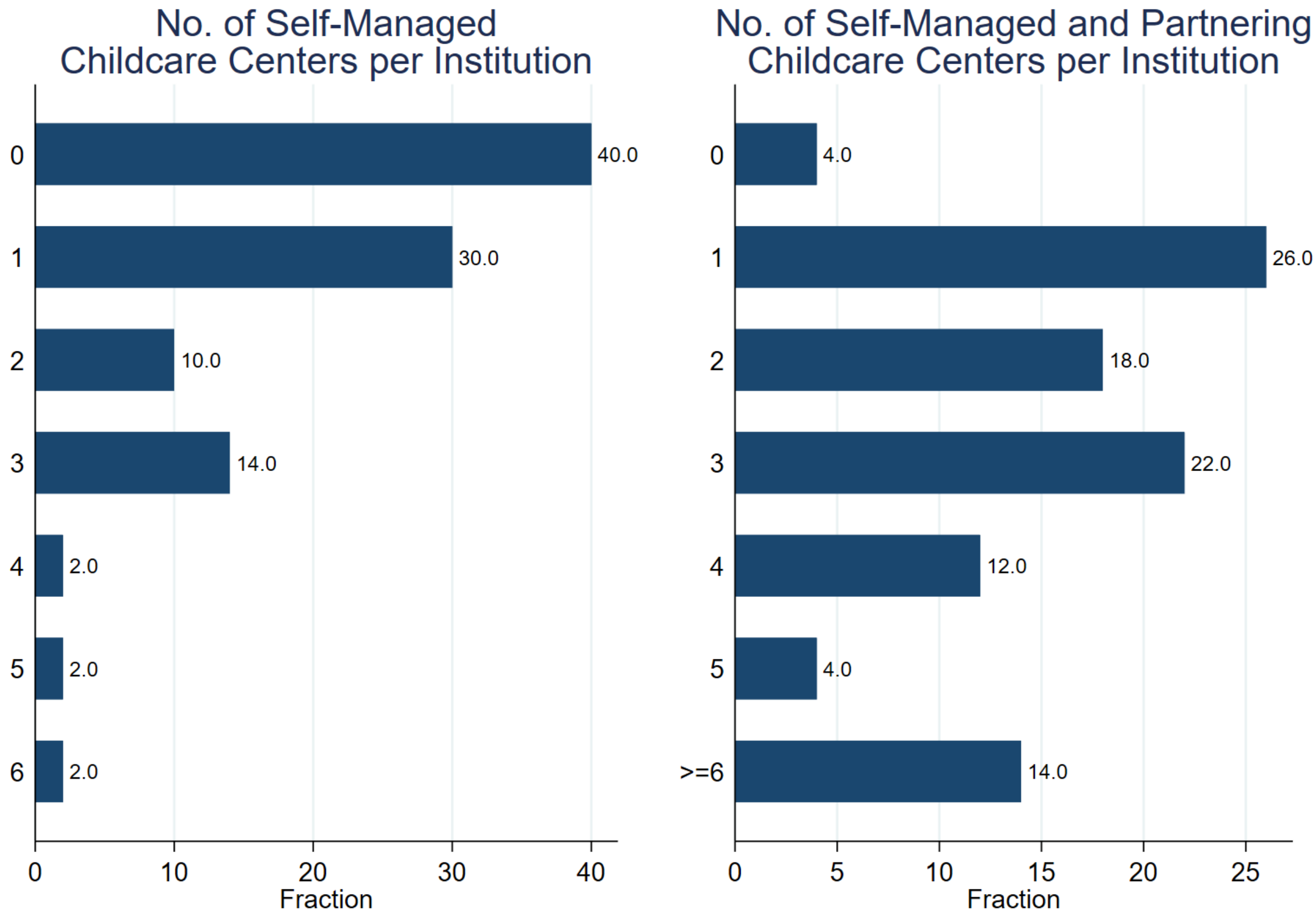

**Fig. S1. Number of childcare centers at the top 50 US research institutions ranked by R&D expenditures.** The figure displays histograms showing the number of childcare centers at the top 50 U.S. institutions, ranked by research and development expenditures. The list of institutions is obtained from the National Science Foundation. The 2022 ranking is used. Percentages are shown next to the bars. The figure on the left shows the histogram of self-managed childcare centers per institution. These centers are typically managed by universities' human resources departments. The figure on the right shows the histogram of self-managed and partnering childcare centers per institution. Partnering centers are those managed by either independent childcare providers or childcare chains. These centers may provide university affiliates (i.e., faculty, staff, students) childcare benefits such as priority enrollment and tuition discounts. Two institutions have neither self-managed nor affiliated centers providing regular care; one has an affiliated center for backup care only. Data are obtained in summer 2024.

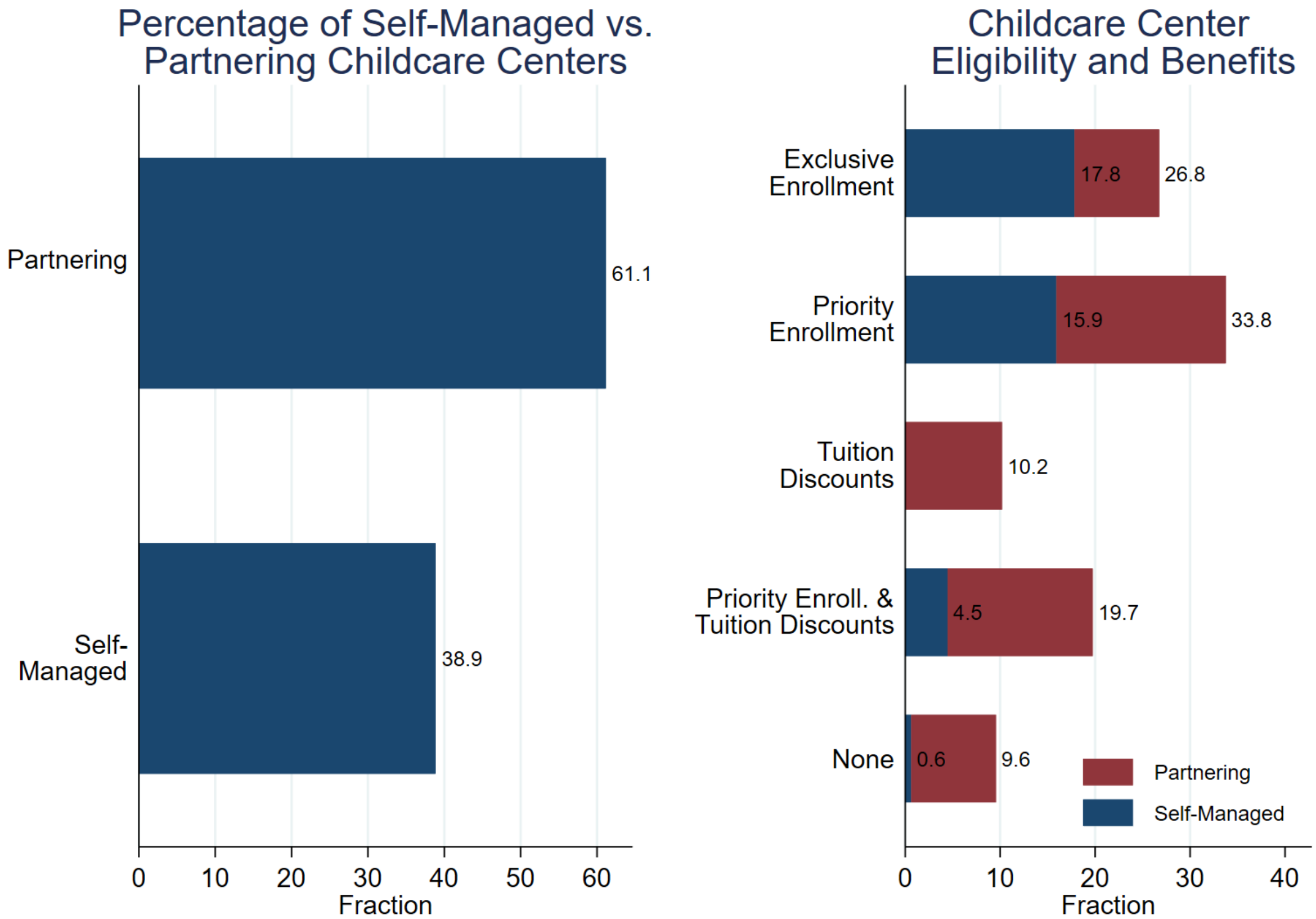

**Fig. S2. Childcare center types, enrollment eligibility, and benefits at the top 50 US research institutions.** Figure on the left displays the percentages of childcare centers managed by the institution itself (38.9%) and by partnering providers (61.1%). Sample is all childcare centers at top 50 US research institutions ranked by R&D expenditures in 2022. $N = 157$ centers. The figure on the right displays the percentages of centers offering exclusive enrollment (26.8%), priority enrollment (33.8%), tuition discounts (10.2%), and priority enrollment and tuition discounts (19.7%) to university affiliates including faculty, staff, and students. **Exclusive enrollment** means only university affiliates (faculty, staff, and students) are eligible to enroll their children in the center. **Priority enrollment** means enrollment is open to the community; however, university affiliates receive waiting list priority (or dedicated space) for enrollment. **Tuition discounts** means enrollment is open to the community; however, university affiliates receive tuition discounts for their enrolled children. Typical discount levels are 5%, 10%, and 15% off. Tuition discounts may be offered simultaneously with priority enrollment. Data are obtained in summer 2024.

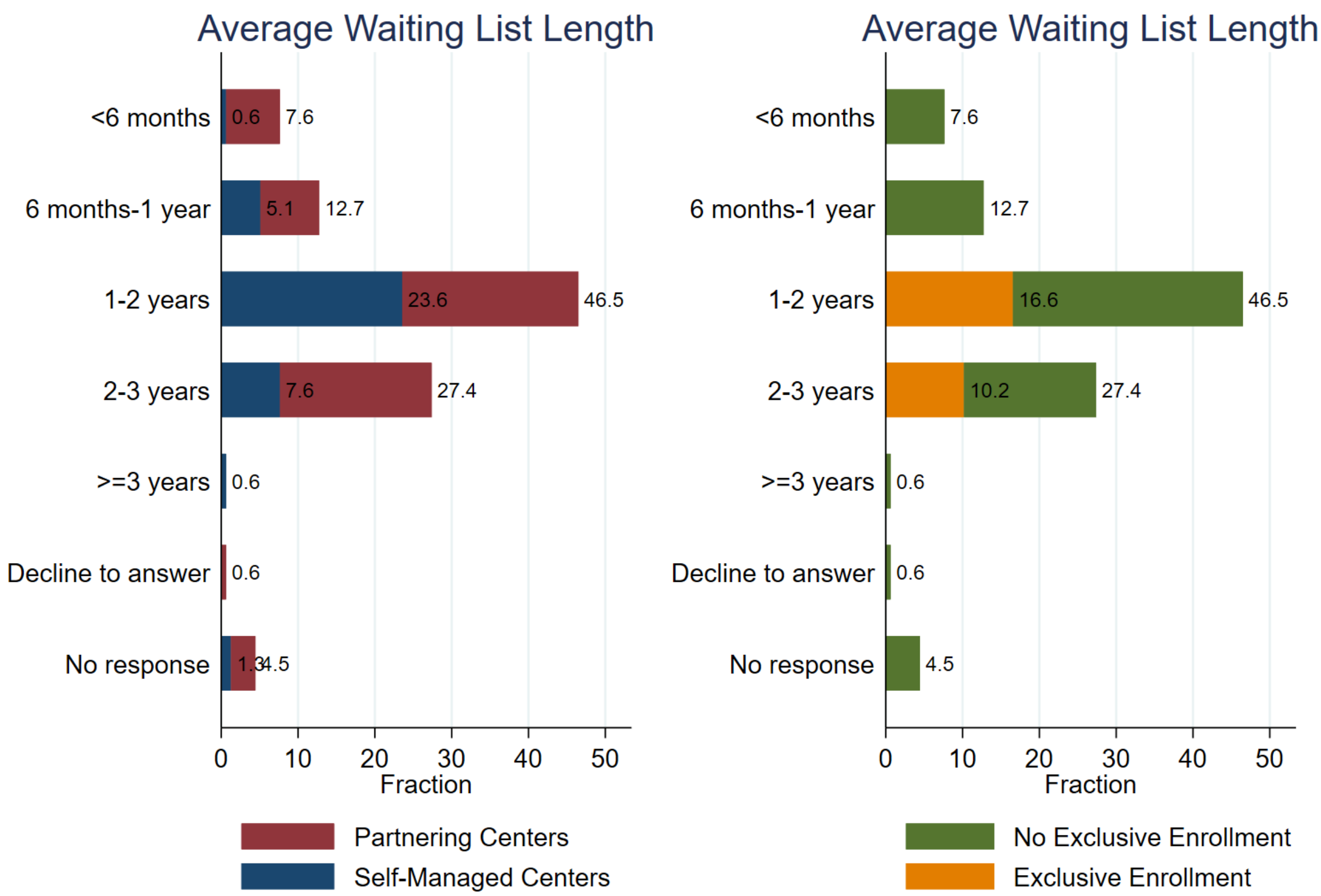

**Fig. S3. Average waiting-list length by childcare center type and enrollment eligibility.** The figure on the left displays center-level average waiting list lengths for self-managed and partnering childcare centers at top 50 US research institutions. $N = 157$ centers. Percentages use all 157 centers and are displayed next to the bars. The largest share of centers (46.5%) have an average waiting list length of 1 (incl.) to 2 years, followed by 27.4% of centers with an average waiting length of 2 (incl.) to 3 years, and 12.7% of centers with an average waiting length of 6 months (incl.) to 1 year. A further 7.6% of centers have waiting lists shorter than 6 months. One center has a waiting list of at least 3 years. Moreover, 4.5% of centers did not respond to the inquiry, and one center declined to respond. Assigning 3, 9, 18, and 30 months to the categories below 6 months, 6 months to 1 year, 1 to 2 years, and 2 to 3 years, respectively, and 42 months to waits of ≥ 3 years, we obtain a mean waiting time of 19.2 months among the 149 centers with reported waiting times ($SD = 8.6$ months). The seven centers with no response and the one center that declined to respond are excluded from these summaries. We do not observe a statistically significant difference between self-managed and partnering centers. The figure on the right displays center-level average waiting list lengths for centers with and without exclusive enrollment policy. Centers with exclusive enrollment for university affiliates have a longer mean waiting time than those without exclusive enrollment policies: 22.6 months (42 centers) versus 17.9 months (107 centers) ($P < 0.05$). Waiting list length data are obtained through center webpage searches (12%), email inquiries (22%) and phone (66%) inquiries in summer 2024.

**Table S1. Sample summary statistics: familial, educational, employment, and regional characteristics of science, engineering, and health doctoral recipients who graduated in 1997–2014 and worked at US universities and research institutions.**

| | Men | Women | Total | *P* |
|---|---|---|---|---|
| Childcare desert (%) | 57.5 | 58.0 | 57.7 | 0.592 |
| Tenured (%) | 38.4 | 32.6 | 35.5 | < 0.001 |
| *Individual characteristics* | | | | |
| Has child aged 0–5 (%) | 46.3 | 39.4 | 42.8 | < 0.001 |
| Race/ethnicity | | | | |
| Asian (%) | 24.6 | 20.3 | 22.4 | < 0.001 |
| Black (%) | 7.5 | 9.3 | 8.4 | < 0.001 |
| Hispanic (%) | 9.1 | 9.5 | 9.3 | 0.506 |
| Native American (%) | 1.8 | 1.9 | 1.9 | 0.537 |
| US citizen (%) | 66.0 | 75.1 | 70.6 | < 0.001 |
| Married or partnered (%) | 76.1 | 72.7 | 74.4 | < 0.001 |
| Spouse working | | | | < 0.001 |
| Full time (%) | 52.3 | 83.5 | 67.7 | |
| Part time (%) | 15.1 | 7.5 | 11.4 | |
| Not working (%) | 32.6 | 9.0 | 20.9 | |
| Age at PhD completion | 32.08 | 32.85 | 32.47 | < 0.001 |
| | (5.84) | (7.13) | (6.53) | |
| PhD at highest research activity inst. (%) | 76.5 | 74.1 | 75.3 | 0.001 |
| PhD field of study | | | | < 0.001 |
| Agricultural and natural resources | 5.4 | 5.2 | 5.3 | |
| Biological and biomedical sciences | 24.7 | 27.4 | 26.0 | |
| Engineering | 15.9 | 6.1 | 11.0 | |
| Health sciences | 5.3 | 10.4 | 7.9 | |
| Mathematics and computer sciences | 9.7 | 5.3 | 7.5 | |
| Physical sciences | 17.4 | 11.7 | 14.5 | |
| Psychology, social sciences, humanities, and other | 21.5 | 33.9 | 27.8 | |

**Table S1. Sample summary statistics: familial, educational, employment, and regional characteristics of science, engineering, and health doctoral recipients who graduated in 1997–2014 and worked at US universities and research institutions (continued).**

| | Men | Women | Total | *P* |
|---|---|---|---|---|
| *Employer characteristics* | | | | |
| Highest research activity institution (%) | 55.2 | 51.1 | 53.1 | < 0.001 |
| Faculty count (log scaled) | | | | |
| Total | 6.96 (1.09) | 6.89 (1.13) | 6.93 (1.11) | < 0.001 |
| Full rank | 5.77 (1.19) | 5.69 (1.24) | 5.73 (1.22) | < 0.001 |
| Associate rank | 5.45 (1.00) | 5.41 (1.04) | 5.43 (1.02) | 0.052 |
| Assistant rank | 5.56 (1.03) | 5.54 (1.08) | 5.55 (1.06) | 0.192 |
| Percentage of women faculty | | | | |
| Total | 38.0 | 40.2 | 39.1 | < 0.001 |
| Full rank | 23.8 | 25.8 | 24.8 | < 0.001 |
| Associate rank | 37.9 | 40.0 | 39.0 | < 0.001 |
| Assistant rank | 44.5 | 46.8 | 45.7 | < 0.001 |
| Student count (log scaled) | | | | |
| Undergraduate | 9.05 (1.93) | 8.96 (2.09) | 9.00 (2.02) | 0.008 |
| Graduate | 8.07 (2.11) | 8.03 (2.12) | 8.05 (2.12) | 0.366 |
| *Regional labor market characteristics* | | | | |
| Population density (log) | 5.76 (1.30) | 5.88 (1.26) | 5.82 (1.28) | < 0.001 |
| Business density (log) | 2.27 (1.13) | 2.36 (1.11) | 2.31 (1.12) | < 0.001 |
| Employment density (log) | 4.83 (1.39) | 4.95 (1.34) | 4.89 (1.37) | < 0.001 |
| Employment annual payroll density (log) | 8.57 (1.64) | 8.73 (1.58) | 8.65 (1.61) | < 0.001 |
| Median household income (log) | 10.90 (0.24) | 10.93 (0.24) | 10.91 (0.24) | < 0.001 |
| Unemployment rate | 5.11 (1.60) | 5.07 (1.51) | 5.09 (1.56) | 0.099 |
| Age 0–5 population density (log) | 3.07 (1.19) | 3.17 (1.17) | 3.12 (1.18) | < 0.001 |
| Democratic presidential vote share (%) | 54.2 (11.5) | 54.8 (11.3) | 54.5 (11.4) | 0.011 |
| Observations | 6,243 | 6,386 | 12,629 | |

*Notes.* All statistics in this table are unweighted. Data sources: (1) childcare deserts: Infogroup US Historical Business Data; (2) tenure and individual characteristics: Survey of Doctorate Recipients (SDR) and Doctorate Records File (DRF); (3) employer characteristics: Integrated Postsecondary Education Data System (IPEDS); (4) regional population and economic characteristics: the Census Bureau; (5) county Democratic presidential vote shares: Presidential Elections in US Counties (1868–2020). Childcare deserts are census tracts with a childcare supply-to-demand ratio below one third. Having a child aged 0–5 is measured within the first six years after doctoral graduation. Tenure is observed at least seven years after doctoral graduation. Highest research activity institutions are defined by the Carnegie Classification of Institutions of Higher Education. Regional population and economic characteristics are measured at the commuting-zone level; Democratic presidential vote shares are measured at the county level. See Methods for sample construction and variable definitions. Standard deviations for continuous variables are in parentheses. Results are consistent when counties are used instead of commuting zones.

**Table S2. Childcare-desert exposure and parenthood timing are similar across gender and parenthood groups.**

| Comparison | Group order (A/B) | $n_A/n_B$ | A/B | $P$ |
|---|---|---|---|---|
| **A. Childcare-desert exposure (%)** | | | | |
| Gender | Men/ | 6,243/ | 57.5/ | 0.592 |
| | Women | 6,386 | 58.0 | |
| Parenthood, all | No child/ | 7,537/ | 57.4/ | 0.446 |
| | Has child | 5,092† | 58.1 | |
| Parenthood, women | No child/ | 4,038/ | 57.6/ | 0.414 |
| | Has child | 2,348 | 58.6 | |
| Parenthood, men | No child/ | 3,499/ | 57.3/ | 0.746 |
| | Has child | 2,744 | 57.7 | |
| **B. Timing of first observed parenthood [mean (SD), years]** | | | | |
| Desert status, all | Non-desert/ | 2,508/ | 4.32 (2.61)/ | 0.811 |
| | Desert | 3,456 | 4.30 (2.55) | |
| Desert status, women | Non-desert/ | 1,139/ | 4.38 (2.48)/ | 0.494 |
| | Desert | 1,596 | 4.32 (2.37) | |
| Desert status, men | Non-desert/ | 1,369/ | 4.26 (2.71)/ | 0.809 |
| | Desert | 1,860 | 4.29 (2.70) | |

*Notes.* All statistics in this table are unweighted. Statistics separated by a slash correspond to groups A and B in the order shown in the second column, Group order (A/B). Column $n_A/n_B$ shows sample sizes. Panel A reports the percentage exposed to a childcare desert. Panel B reports means and standard deviations in parentheses. †A total of 5,408 (42.8% out of N=12,629; table S1) scientists had young children below age 6 within the six years post doctoral graduation. Among them, 5,092 scientists had children while they were employed at universities or research institutions. The remaining 316 scientists had children while they were not working at universities or research intuitions.

**Table S3. Associations of gender, parenthood, and workplace childcare deserts with scientists' tenure attainment.**

| Sample | All scientists | | Parents | | Nonparents | |
|---|---|---|---|---|---|---|
| | (1) | (2) | (3) | (4) | (5) | (6) |
| Women | $-0.045^{***}$ | $-0.026^{**}$ | $-0.060^{***}$ | $-0.020$ | $-0.025$ | $-0.016$ |
| | (0.010) | (0.012) | (0.015) | (0.024) | (0.016) | (0.025) |
| Parenthood | | $-0.002$ | | | | |
| | | (0.017) | | | | |
| Women × parenthood | | $-0.049^{***}$ | | | | |
| | | (0.018) | | | | |
| Childcare desert | | | | 0.029 | | $-0.002$ |
| | | | | (0.023) | | (0.023) |
| Women × childcare desert | | | | $-0.069^{**}$ | | $-0.014$ |
| | | | | (0.032) | | (0.028) |
| Observations | 12,150$^{\dagger}$ | 12,150 | 5,023 | 5,011$^{\ddagger}$ | 5,897 | 5,879 |
| Adjusted $R^2$ | 0.18 | 0.18 | 0.22 | 0.22 | 0.17 | 0.17 |

*Notes.* The model is ordinary least squares (OLS) regression. The outcome variable is the binary indicator of whether an individual has attained tenure promotion observed at least seven years after doctoral graduation. All models include the same controls as in Fig. 2C. Standard errors are in parentheses and are clustered at the state level. $^{\dagger}$ 479 (3.9%) observations are dropped due to missing values for employer characteristics: some institutions did not report faculty counts or enrollments in certain years. $^{\ddagger}$ 12 (0.2%) observations are dropped due to missing child care availability data caused by missing or inaccurate facility location information.

$^{***}P < 0.01$; $^{**}P < 0.05$; $^{*}P < 0.10$.

**Table S4. Childcare deserts in focal and neighboring tracts and tenure promotion.**

| Desert specification | Separate focal and neighboring tracts | |
|---|---|---|
| Sample | Parents | Nonparents |
| | (1) | (2) |
| Women | 1.016 | 0.931 |
| | (0.093) | (0.100) |
| Childcare desert: focal tract | 1.030 | 0.979 |
| | (0.087) | (0.085) |
| Women × childcare desert: focal tract | $0.750^{**}$ | 1.020 |
| | (0.085) | (0.126) |
| Childcare desert: neighboring tracts | 1.129 | 1.023 |
| | (0.099) | (0.077) |
| Women × childcare desert: neighboring tracts | $0.821^{**}$ | 0.967 |
| | (0.067) | (0.109) |
| Observations | 21,674 | 24,831 |
| Log pseudolikelihood | −472,910.40 | −485,208.02 |

*Notes.* Each column reports a survival model of tenure promotion that jointly includes separate childcare-desert indicators for the focal census tract and neighboring census tracts, together with their interactions with women. All models include the same controls as in Fig. 2C. Standard errors are clustered at the state level. Column (1) uses the parent sample; column (2) uses the nonparent sample. Hazard ratios are reported, with standard errors in parentheses. Hazard ratios below 1 indicate a lower instantaneous rate of tenure promotion. Observation counts refer to follow-up records.

$^{***}P < 0.01$; $^{**}P < 0.05$; $^{*}P < 0.10$.